\documentclass[aos,MSNbibl,nameyear,dvips]{arximspdf}
\usepackage[ruled,vlined]{algorithm2e}
\usepackage{graphicx}

\doi{10.1214/13-AOS1119} 
\volume{41}
\issue{3}
\pubyear{2013}
\firstpage{1516}
\lastpage{1541}

\makeatletter

\newcommand{\underset}[2]{\,\mathop{#2}_{#1}\,}

\newcommand{\Ed}{\operatorname{E}}

\newcommand{\cD}{\mathcal{D}}
\newcommand{\cK}{\mathcal{K}}
\newcommand{\cC}{\mathcal{C}}
\newcommand{\cF}{\mathcal{F}}
\newcommand{\cS}{\mathcal{S}}

\renewcommand{\P}{\mathbb{P}}
\newcommand{\E}{\mathbb{E}}
\newcommand{\R}{\mathbb{R}}

\newcommand{\oI}{\overline{I}}
\newcommand{\ind}{\mathbb{I}}
\newcommand{\measg}{\mathfrak{M}_1}

\newcommand{\argmax}{\mathop{\arg\max}}
\newcommand{\KL}{\operatorname{KL}}

\renewcommand{\d}{\mathrm{d}}
\newcommand{\defeq}{\stackrel{\mathrm{def}}{=}}
\newcommand{\eqdef}{\defeq}
\newcommand{\Supp}{\operatorname{Supp}}
\newcommand{\Vard}{\operatorname{Var}}

\newcommand{\ber}{\mathrm{BER}}
\newcommand{\qua}{\mathrm{QUAD}}

\renewcommand{\epsilon}{\varepsilon}

\newproclaim{example}{Example}
\newtheorem{theorem}{Theorem}
\newtheorem{corollary}{Corollary}
\newtheorem{lemma}{Lemma}
\newtheorem{proposition}{Proposition}
\newtheorem{toprove}{Fact to be proven}

\makeatother

\begin{document}
\begin{frontmatter}

\title{Kullback--Leibler upper confidence bounds for optimal sequential allocation}
\runtitle{Kullback--Leibler upper confidence bounds}

\begin{aug}
\author[A]{\fnms{Olivier} \snm{Capp\'e}\ead[label=e1]{cappe@telecom-paristech.fr}},
\author[B]{\fnms{Aur\'elien} \snm{Garivier}\corref{}\ead[label=e2]{aurelien.garivier@math.univ-toulouse.fr}},
\author[C]{\fnms{Odalric-Ambrym}~\snm{Maillard}\thanksref{t1}\ead[label=e3]{odalricambrym.maillard@gmail.com}},
\author[D]{\fnms{R\'emi}~\snm{Munos}\thanksref{t1}\ead[label=e4]{remi.munos@inria.fr}}
\and
\author[E]{\fnms{Gilles} \snm{Stoltz}\thanksref{t1}\ead[label=e5]{gilles.stoltz@ens.fr}}
\runauthor{O. Capp\'e et al.}
\affiliation{Telecom ParisTech, CNRS,
Universit\'{e} Paul Sabatier, University of Leoben, INRIA, Ecole Normale
Sup\'{e}rieure, HEC Paris}
\address[A]{O. Capp\'e\\
LTCI, Telecom ParisTech, CNRS \\
46 rue Barrault\\
75013 Paris\\
France \\
\printead{e1}}
\address[B]{A. Garivier\\
IMT\\
Universit\'e Paul Sabatier \\
118 route de Narbonne\\
31000 Toulouse\\
France \\
\printead{e2}}
\address[C]{O.-A. Maillard\\
University of Leoben \\
8700 Leoben\\
Austria \\
\printead{e3}}
\address[D]{R. Munos\\
INRIA Lille, SequeL Project\hspace*{57.45pt} \\
40 avenue Halley\\
59650 Villeneuve d'Ascq\\
France \\
\printead{e4}}
\address[E]{G. Stoltz\\
Ecole Normale Sup\'erieure\\
CNRS, INRIA \\
45 rue d'Ulm, 75005 Paris\\
France \\
and \\
HEC Paris, CNRS \\
1 rue de la Lib\'eration\\
78351 Jouy-en-Josas\\
France \\
\printead{e5}} 
\end{aug}

\thankstext{t1}{Supported by the French National Research Agency
(ANR) under Grant EXPLO/RA (``Exploration--exploitation for efficient
resource allocation''), the PASCAL2 Network of Excellence under EC
Grant 506778 and the ECs Seventh Framework Programme
(FP7/2007--2013) under Grant agreement no. 270327.}

\received{\smonth{10} \syear{2012}}
\revised{\smonth{3} \syear{2013}}

%
\begin{abstract}
We consider optimal sequential allocation in the context of the
so-called stochastic multi-armed bandit model. We describe a generic
index policy, in the sense of Gittins [\textit{J. R. Stat. Soc. Ser. B
Stat. Methodol.} \textbf{41} (1979) 148--177], based on upper
confidence bounds of the arm payoffs computed using the
Kullback--Leibler divergence. We consider two classes of distributions
for which instances of this general idea are analyzed: the
\texttt{kl-UCB} algorithm is designed for one-parameter
exponential families and the empirical \texttt{KL-UCB} algorithm for bounded and
finitely supported distributions. Our main contribution is a unified
finite-time analysis of the regret of these algorithms that
asymptotically matches the lower bounds of Lai and Robbins
[\textit{Adv. in Appl. Math.} \textbf{6} (1985) 4--22] and Burnetas and
Katehakis [\textit{Adv. in Appl. Math.} \textbf{17} (1996) 122--142],
respectively. We also investigate the behavior of these algorithms when
used with general bounded rewards, showing in particular that they
provide significant improvements over the state-of-the-art.
\end{abstract}

%

\end{frontmatter}

\section{Introduction}\label{intro}

This paper is about optimal sequential allocation in unknown random
environments. More precisely, we consider the setting known under the
conventional, if not very explicit, name of (stochastic) \emph{multi-armed
bandit}, in reference to the 19th century gambling game. In the
multi-armed bandit
model, the emphasis is put on focusing as quickly as possible on the best
available option(s) rather than on estimating precisely the efficiency
of each
option. These options are referred to as arms, and each of them is
associated with a distribution;
arms are indexed by $a$ and associated distributions are denoted by
$\nu_a$.

The archetypal example occurs in clinical trials where the options (or
arms) correspond to available treatments whose efficiencies are unknown
a priori, and patients arrive sequentially; the action consists of
prescribing a particular treatment to the patient, and the observation
corresponds (e.g.) to the success or failure of the treatment. The goal
is clearly here to achieve as many successes as possible. A~strategy
for doing so is said to be \emph{anytime} if it does not require to
know in advance the number of patients that will participate to the
experiment. Although the term multi-armed bandit was probably coined in
the late 1960s [\citet{gittins-1979}], the origin of the problem
can be traced back to fundamental questions about optimal stopping
policies in the context of clinical trials [see Thompson
(\citeyear{Tho33,Tho35})] raised since the 1930s; see also
\citet{wald-1945,ro52}.

In his celebrated work, \citet{gittins-1979} considered the
\emph{Bayesian-optimal} solution to the discounted infinite-horizon multi-armed
bandit problem. Gittins first showed that the Bayesian optimal policy
could be
determined by dynamic programming in an extended Markov decision
process. The
second key element is the fact that the optimal policy search can be factored
into a set of simpler computations to determine \emph{indices} that fully
characterize each arm given the current history of the game
[\citet{gittins-1979,whittle-1980,weber-1992}]. The optimal
policy is then an
\emph{index policy} in the sense that at each time round, the (or an)
arm with
highest index is selected. Hence, index policies only differ in the way the
indices are computed.

From a practical perspective, however, the use of Gittins indices is
limited to
specific arm distributions and is computationally challenging
[\citet{gittins2011multi}].
In the 1980s, pioneering works by \citet{LaiRo85},
\citet{chang-lai-1987}, Burnetas and Katehakis
(\citeyear{Burnetas96,BurnetasKatehakis97,burnetas-katehakis-2003}) suggested that
Gittins indices can be approximated by
quantities that can be interpreted as upper bounds of confidence intervals.
\citet{Agrawal95} formally introduced and provided an asymptotic
analysis for generic classes of index policies termed \texttt{UCB} (for
Upper Confidence Bounds).
For general bounded reward distributions, \citet{Auer02} provided
a finite time analysis for a particular variant of \texttt{UCB} based on
Hoeffding's inequality; see also \citet
{bubeckcesa-bianchi2012-ftml} for a recent survey of bandit models
and variants.

There are, however, significant differences between the algorithms and
results of
\citet{gittins-1979} and \citet{Auer02}. First, \texttt{UCB} is
an anytime algorithm that does not rely on the
use of a discount factor or even on the knowledge of the horizon of the
problem. More significantly, the Bayesian perspective is absent, and
\texttt{UCB} is
analyzed in terms of its frequentist (distribution-dependent or
distribution-free)
performance, by exhibiting \emph{finite-time}, nonasymptotic
bounds on its expected regret. The expected regret of an algorithm---a
quantity to be formally defined in
Section~\ref{secnot}---corresponds to the difference, in expectation,
between the rewards that would have been
gained by only pulling a best arm and the rewards actually gained.

\texttt{UCB} is a very robust algorithm that is suited to all problems with
bounded stochastic rewards
and has strong performance guarantees, including distri\-bu\-tion-free ones.
However, a closer examination of the arguments in the proof reveals
that the form of the upper confidence bounds used in \texttt{UCB} is a direct
consequence of the use of Hoeffding's inequality and
significantly differs from the approximate form of Gittins indices suggested
by \citet{LaiRo85} or \citet{Burnetas96}. Furthermore, the
frequentist
asymptotic lower bounds for the regret obtained by these authors also
suggest that the behavior of \texttt{UCB} can be far from optimal. Indeed,
under suitable
conditions on the model $\cD$ (the class of possible distributions
associated with each arm),
any policy that is ``admissible'' [i.e., not grossly under-performing;
see \citet{LaiRo85} for details] must satisfy
the following asymptotic inequality on its expected regret $\E[R_T]$
at round $T$:
%
\begin{equation}
\label{eqbinfRegret}\liminf_{T\to\infty} \frac{\E[R_T]}{\log
(T)} \geq\sum
_{a:\mu_a<\mu^\star} \frac{\mu^\star-\mu_a}{\cK_{\inf} (\nu
_a,\mu^\star)},
\end{equation}
where $\mu_a$ denotes the expectation of the
distribution $\nu_a$ of arm $a$, while $\mu^\star$ is the maximal
expectation among all arms. The
quantity
%
\begin{equation}
\label{eqKinfdef} \cK_{\inf} (\nu,\mu) = \inf \bigl\{ \KL \bigl(\nu,
\nu' \bigr)\dvtx \nu' \in\cD\mbox{ and } \Ed \bigl(
\nu' \bigr) > \mu \bigr\},
\end{equation}
which measures the difficulty of the problem, is the minimal Kullback--Leibler
divergence between the arm distribution $\nu$ and distributions in the model
$\cD$ that have expectations larger than $\mu$. By comparison, the bound
obtained in \citet{Auer02} for \texttt{UCB} is of the form
\[
\E[R_T] \leq C \biggl( \sum_{a:\mu_a<\mu^\star}
\frac{1}{\mu^\star-\mu_a} \biggr) \log(T)+ o \bigl(\log(T) \bigr)
\]
for some numerical constant $C$, for example, $C=8$; we provide a
refinement of the result of
\citet{Auer02} as Corollary~\ref{corUCB}, below. These two
results coincide
as to the logarithmic rate of the expected regret, but the
(distribution-dependent) constants differ, sometimes
significantly. Based on this observation, Honda and
Takemura (\citeyear{HondaTakemura10DMED,HondaTakemura11ML}) proposed
an algorithm, called \texttt{DMED}, that is not an index policy but
was shown to improve over \texttt{UCB} in some situations.
They later showed that this algorithm could also accommodate the case
of semi-bounded rewards; see \citet{Honda2012}.

Building on similar ideas, we show in this paper that for a large class
of problems there does
exist a generic index policy---following the insights of \citet
{LaiRo85}, \citet{Agrawal95}
and \citet{Burnetas96}---that guarantees a bound on
the expected regret of the form
\[
\E[R_T] \leq\sum_{a:\mu_a<\mu^\star} \biggl(
\frac{\mu^\star-\mu_a}{\cK_{\inf} (\nu_a,\mu^\star)} \biggr) \log(T)+ o \bigl(\log(T) \bigr),
\]
and which is thus asymptotically optimal.\setcounter{footnote}{1}\footnote{Minimax optimality
is another, distribution free, notion of optimality that has also been
studied in the bandit setting [\citet
{bubeckcesa-bianchi2012-ftml}]. In this paper, we focus on
problem-dependent optimality.} Interestingly, the index used in this
algorithm can be interpreted as the upper bound of a confidence region
for the
expectation constructed using an empirical likelihood principle
[\citet{owen2001empirical}].

We describe the implementation of this algorithm and analyze its
performance in
two practically important cases where the lower bound of (\ref{eqbinfRegret})
was shown to hold [\citet{LaiRo85,Burnetas96}]---namely, for
one-parameter canonical
exponential families of distributions (Section~\ref{secklUCB}), in
which case
the algorithm is referred to as \texttt{kl-UCB}, and for
finitely supported distributions
(Section~\ref{secKLemp}), where the algorithm is called empirical
\texttt{KL-UCB}. Determining
the empirical \texttt{KL-UCB} index requires solving a convex program
(maximizing a linear
function on the probability simplex under Kullback--Leibler
constraints) for
which we provide in the supplemental article [\citet{extended},
Appendix C.1]
a simple algorithm inspired by \citet{FilippiCappeGarivier10}.

The analysis presented here greatly improves over the preliminary
results presented, on the one hand by \citet{GaCa11}, and on the
other hand by \citet{MaillardMunosStoltz11klucb}; more precisely,
the improvements lie in the greater generality of the analysis and by
the more precise evaluation of the remainder terms in the regret
bounds. We believe that the result obtained in this paper for
\texttt{kl-UCB} (Theorem~\ref{thmklUCB}) is not improvable. For
empirical \texttt{KL-UCB} the bounding of the remainder term could be
improved upon obtaining a sharper version of the contraction lemma for
$\cK_{\inf}$ [Lemma 6 in the supplemental article,
\citet{extended}]. The proofs rely on results of independent
statistical interest: nonasymptotic bounds on the level of sequential
confidence intervals for the expectation of independent, identically
distributed variables, (1) in canonical exponential families
(equation~(\ref{eqinegDevIndepInterest}); see also Lemma 11 in the
supplemental article [\citet{extended}]) and (2) using the empirical
likelihood method for bounded variables (Proposition~\ref {propcovEL}).

For general bounded distributions, we further make three important
observations. First, the particular instance of the
\texttt{kl-UCB} algorithm based
on the Kullback--Leibler divergence\vadjust{\goodbreak} between normal distributions is the
\texttt{UCB}
algorithm, which allows us to provide an improved optimal finite-time
analysis of
its performance (Corollary~\ref{corUCB}). Next, the
\texttt{kl-UCB} algorithm, when used with
the Kullback--Leibler divergence between Bernoulli distributions, obtains
a strictly better performance than \texttt{UCB}, for any bounded
distribution (Corollary~\ref{thborneRegretklUCB}). Finally,
although a complete analysis of the empirical \texttt{KL-UCB} algorithm is
subject to further investigations, we show here
that the empirical \texttt{KL-UCB} index has a guaranteed coverage probability
for general bounded distributions, in the sense that, at any step, it
exceeds the true expectation with large probability (Proposition \ref
{propcovEL}). We
provide some empirical evidence that empirical \texttt{KL-UCB} also performs
well for
general bounded distributions and illustrate the tradeoffs arising when using
the two algorithms, in particular for short horizons.

\subsection*{Outline}
The paper is organized as follows. Section~\ref{secnot} introduces the
necessary notations and defines the notion of regret. Section~\ref{secKLUCB}
presents the generic form of the \texttt{KL-UCB} algorithm and provides the
main steps
for its analysis, leaving two facts to be proven under each specific
instantiation
of the algorithm. The \texttt{kl-UCB} algorithm in the case of
one-dimensional exponential families is considered in
Section~\ref{secklUCB}, and the empirical \texttt{KL-UCB} algorithm for
bounded and
finitely supported distributions is presented in
Section~\ref{secKLemp}. Finally, the behavior of these algorithms in
the case
of general bounded distributions is investigated in Section~\ref{secgeneral};
and numerical experiments comparing \texttt{kl-UCB} and
empirical \texttt{KL-UCB} to their
competitors are reported in Section~\ref{secexpnum}. Proofs are provided
in the supplemental article [\citet{extended}].

\section{Setup and notation}
\label{secnot}

We consider a bandit problem with finitely many arms indexed by $a \in
\{ 1,\ldots, K \}$, with $K \geq2$,
each associated with an (unknown) probability distribution $\nu_a$
over $\R$.
We assume, however, that a model $\cD$ is known: a family of
probability distributions such that
$\nu_a \in\cD$ for all arms $a$.

The game is sequential and goes as follows:
at each round $t \geq1$, the player picks an arm $A_t$ (based on the
information gained in the past)
and receives a stochastic payoff $Y_t$ drawn independently
at random according to the distribution $\nu_{A_t}$. He only gets to
see the payoff $Y_t$.

\subsection{Assessment of the quality of a strategy via its expected regret}
For each arm $a \in\{ 1,\ldots,K \}$, we denote by $\mu_a$ the
expectation of its associated distribution~$\nu_a$,
and we let $a^\star$ be any optimal arm, that is,
\[
a^{\star} \in\argmax_{a \in\{ 1,\ldots,K \} } \mu_a.
\]
We write $\mu^\star$ as a short-hand notation for the largest
expectation $\mu_{a^{\star}}$
and denote the gap of the expected payoff $\mu_a$ of an arm $a$ to
$\mu^\star$ as
$\Delta_a = \mu^\star- \mu_a$.
In addition, the number of times each arm $a$ is pulled between the
rounds $1$ and $T$
is referred to as $N_a(T)$,
\[
N_a(T) \eqdef\sum_{t=1}^T
\ind_{ \{ A_t = a \} }.
\]

The quality of a strategy will be evaluated through the standard notion
of expected regret, which
we define formally now. The expected regret (or simply, regret) at
round $T \geq1$ is defined as
%
\begin{equation}
\label{eqdefregr} R_T \eqdef\E \Biggl[ T \mu^{\star} - \sum
_{t=1}^T Y_t \Biggr] = \E
\Biggl[ T \mu^{\star} - \sum_{t=1}^T
\mu_{A_t} \Biggr] = \sum_{a = 1}^K
\Delta_a \E \bigl[ N_a(T) \bigr],
\end{equation}
where we used the tower rule for the first equality. Note that the
expectation is with respect to the random draws of the $Y_t$ according
to the
$\nu_{A_t}$ and also to the possible auxiliary randomizations that the
decision-making strategy
is resorting to.

The regret measures the cumulative loss resulting from pulling
suboptimal arms, and thus quantifies
the amount of exploration required by an algorithm in order to find a
best arm, since, as (\ref{eqdefregr}) indicates,
the regret scales with the expected number of pulls of suboptimal arms.

\subsection{Empirical distributions}
We will denote them in two related ways, depending on whether random
averages indexed by the global time $t$
or averages of a given number $n$ of pulls of a given arms are
considered. The first series of averages will be
referred to by using a functional notation for the indexing in the
global time: $\widehat{\nu}_a(t)$,
while the second series will be indexed with the local times $n$ in
subscripts: $\widehat{\nu}_{a,n}$.
These two related indexings, functional for global times and random
averages versus
subscript indexes for local times, will be consistent throughout the
paper for all quantities at hand, not only
empirical averages.

More formally, for all arms $a$ and all rounds $t$ such that $N_a(t)
\geq1$,
\[
\widehat{\nu}_a(t) = \frac{1}{N_a(t)} \sum
_{s=1}^t \delta_{Y_s} \ind_{ \{ A_s = a \} },
\]
where $\delta_x$ denotes the Dirac distribution on $x \in\R$.

For averages based on local times we need to introduce stopping times.
To that end, we consider the filtration $(\cF_t)$, where for all $t
\geq1$, the $\sigma$-algebra $\cF_t$ is generated
by $A_1,Y_1,\ldots,A_t,Y_t$. In particular, $A_{t+1}$ and all
$N_{a}(t+1)$ are $\cF_t$-measurable.
For all $n \geq1$, we denote by $\tau_{a,n}$ the round at which $a$
was pulled for the $n$th time; since
\[
\tau_{a,n} = \min \bigl\{ t \geq1\dvtx N_a(t) = n \bigr\},
\]
we see that $ \{ \tau_{a,n} = t \}$ is $\cF_{t-1}$-measurable. That
is, each
random variable $\tau_{a,s}$ is a (predictable) stopping time. Hence,
as shown, for instance, in Chow and Teicher [(\citeyear{ChTe88}), Section 5.3], the random
variables $X_{a,n} = Y_{\tau_{a,n}}$, where $n = 1,2,\ldots\,$,
are independent and identically distributed according to $\nu_a$.
For all arms $a$, we then denote by
\[
\widehat{\nu}_{a,n} = \frac{1}{n} \sum_{k=1}^n
\delta_{X_{a,k}}
\]
the empirical distributions corresponding to local times $n \geq1$.

All in all, we of course have the rewriting
\[
\widehat{\nu}_{a}(t) = \widehat{\nu}_{a,N_a(t)}.
\]

\section{The \texttt{KL-UCB} algorithm}
\label{secKLUCB}

We fix an interval or discrete subset
$\cS\subseteq\R$ and denote by $\mathfrak{M}_1(\cS)$ the set of
all probability
distributions over $\cS$.
For two distributions $\nu,\nu' \in\mathfrak{M}_1(\cS)$, we denote
by $\KL(\nu,\nu')$
their Kullback--Leibler divergence and by $\Ed(\nu)$ and $\Ed(\nu
')$ their expectations.
(This expectation operator is denoted by $\Ed$ while expectations with
respect to underlying
randomizations are referred to as~$\E$.)

\begin{algorithm}[b]
\begin{minipage}{340pt}
\caption{The \texttt{KL-UCB} algorithm (generic form).}
\label{algKLUCB}

\textit{Parameters}: An operator $\Pi_\cD\dvtx \mathfrak{M}_1(\cS)\to
\cD$; a
nondecreasing function\\ $f\dvtx  \mathbb{N} \to\mathbb{R}$\\
\textit{Initialization}: Pull each arm of $\{ 1,\ldots, K \}$ once
\\
\For{$t=K$ to $T-1$,}{

compute for each arm $a$ the quantity
%
\begin{equation}
\label{eqUCBdef} U_a(t) = \sup \biggl\{ \Ed(\nu)\dvtx  \nu\in\cD
\mbox{ and } \KL \bigl( \Pi_\cD \bigl( \widehat{\nu}_a(t) \bigr),
\nu \bigr) \leq\frac{f(t)}{N_a(t)} \biggr\}
\end{equation}
pick an arm $ {A_{t+1} \in\argmax_{a \in\{ 1,\ldots,K \}} U_a(t)}$
}
\end{minipage}
\end{algorithm}

The generic form of the algorithm of interest in this paper is
described as Algorithm~\ref{algKLUCB}.
It relies on two parameters: an operator $\Pi_\cD$ (in spirit, a
projection operator)
that associates with each empirical distribution $\widehat{\nu}_a(t)$ an
element of the model~$\cD$;
and a nondecreasing function $f$, which is typically such that $f(t)
\approx\log(t)$.

At each round $t \geq K$, an upper confidence bound $U_a(t)$ is
associated with the expectation $\mu_a$
of the distribution $\nu_a$ of each arm; an arm $A_{t+1}$ with highest
upper confidence bound is then
played. Note that the algorithm does not need to know the time horizon
$T$ in advance. Furthermore, the \texttt{UCB} algorithm of \citet
{Auer02} may be recovered by replacing\vadjust{\goodbreak} $\KL( \Pi_\cD( \widehat{\nu
}_a(t) ), \nu)$ with a quantity proportional to $ (\Ed(\widehat{\nu
}_a(t))-\Ed(\nu) )^2$; the implications of this observation will be
made more explicit in Section~\ref{secgeneral}.

\subsection{General analysis of performance}
\label{secgalanal}

In Sections~\ref{secklUCB} and~\ref{secKLemp},
we prove non\-asymptotic regret bounds for Algorithm~\ref{algKLUCB} in
two different settings.
These bounds match the asymptotic lower bound (\ref{eqbinfRegret}) in
the sense that, according to (\ref{eqdefregr}), bounding the expected
regret is equivalent to bounding the number of suboptimal draws. We
show that, for any suboptimal arm $a$, we have
\[
\E \bigl[N_a(T) \bigr] \leq\frac{\log(T)}{\cK_{\inf} (\nu_a,\mu
^\star)} \bigl(1+o(1) \bigr),
\]
where the quantity $\cK_{\inf} (\nu_a,\mu^\star)$ was defined in
the \hyperref[intro]{Introduction}.
This result appears as a consequence of nonasymptotic bounds, which
are derived using a common analysis framework detailed in the rest of
this section.

Note that the term $\log(T)/\cK_{\inf} (\nu_a,\mu^\star)$ has an heuristic
interpretation in terms of large deviations, which gives some insight
on the
regret analysis to be presented below. Let $\nu'\in\cD$ be such that
$\Ed(\nu')\geq\mu^\star$, let $X'_1,\ldots,X'_n$ be independent
variables with
distribution $\nu'$ and let $\widehat{\nu}'_n =
(\delta_{X'_1}+\cdots+\delta_{X'_n})/n$. By Sanov's theorem, for a small
neighborhood $\mathcal{V}_a$ of $\nu_a$, the probability that $\widehat
{\nu}'_n$
belongs to $\mathcal{V}_a$ is such that
\[
-\frac{1}{n}\log\P \bigl\{\widehat{\nu}'_n \in
\mathcal{V}_a \bigr\} \underset{n\to\infty} {\longrightarrow} \inf
_{\nu\in
\mathcal{V}_a}\KL \bigl(\nu, \nu' \bigr) \approx\KL \bigl(
\nu_a, \nu' \bigr) \geq\cK_{\inf} \bigl(
\nu_a, \mu^\star \bigr).
\]
In the limit, ignoring the sub-exponential terms, this means that for
$n=\log(T)/\cK_{\inf} (\nu_a,\mu^\star)$, the probability
$\P\{\widehat{\nu}'_n \in\mathcal{V}_a \}$ is smaller than $1/T$.
Hence, $\log(T)/\cK_{\inf} (\nu_a,\mu^\star)$ appears as the
minimal number $n$ of draws ensuring that the probability under any
distribution with expectation at least $\mu^\star$ of the event ``the
empirical distribution of $n$ independent draws belongs to a
neighborhood of $\nu_a$'' is smaller than $1/T$. This event, of course,
has an overwhelming probability under~$\nu_a$. The significance of
$1/T$ as a cutoff value can be understood as follows: if the suboptimal
arm $a$ is chosen along the $T$ draws, then the regret is at most equal
to $(\mu^\star-\mu_a)T$; thus, keeping the probability of this event
under $1/T$ bounds the contribution of this event to the average regret
by a constant. Incidentally, this explains why knowing $\mu^\star$ in
advance does not significantly reduce the number of necessary
suboptimal draws. The analysis that follows shows that the bandit
problem, despite its sequential aspect and the absence of prior
knowledge on the expectation of the arms, is indeed comparable to a
sequence of tests of level $1-1/T$ with null hypothesis $H_0\dvtx
\Ed(\nu')>\mu^\star$ and alternative hypothesis $H_1\dvtx  \nu' = \nu
_a$, for which Stein's lemma [see, e.g., \citet{VdV98}, Theorem
16.12] states that the best error exponent is $\cK_{\inf}
(\nu_a,\mu^\star)$.

Let us now turn to the main lines of the regret proof. By definition of
the algorithm, at rounds $t \geq K$,\vadjust{\goodbreak} one has $A_{t+1} = a$
only if $U_a(t) \geq U_{a^\star}(t)$. Therefore, one has the decomposition
%
\begin{eqnarray}
\label{eqmajoNTa} \{ A_{t+1} = a \} & \subseteq& \bigl\{ \mu
^{\dagger}\geq U_{a^\star}(t) \bigr\} \cup \bigl\{ \mu^{\dagger}<
U_{a^\star}(t) \mbox{ and } A_{t+1} = a \bigr\}
\nonumber\\[-8pt]\\[-8pt]
& \subseteq& \bigl\{ \mu^{\dagger}\geq U_{a^\star}(t) \bigr
\} \cup \bigl\{ \mu^{\dagger}< U_a(t) \mbox{ and }
A_{t+1} = a \bigr\},\nonumber
\end{eqnarray}
where $\mu^{\dagger}$ is a parameter which is taken either equal to
$\mu
^\star$, or slightly smaller when required by technical arguments.
The event $ \{ \mu^{\dagger}< U_a(t) \}$ can be rewritten as
\begin{eqnarray*}
\bigl\{ \mu^{\dagger}< U_a(t) \bigr\} & = & \biggl\{ \exists
\nu' \in\cD\dvtx  \Ed \bigl(\nu' \bigr) > \mu^{\dagger}
\mbox{ and } \KL \bigl( \Pi_\cD \bigl( \widehat{\nu}_a(t)
\bigr), \nu' \bigr) \leq\frac{f(t)}{N_a(t)} \biggr\}
\\
& = &\bigl\{ \widehat{\nu}_a(t) \in\cC_{\mu^{\dagger}, f(t)/N_a(t)} \bigr\}=\{
\widehat{\nu}_{a,N_a(t)} \in\cC_{\mu^{\dagger}, f(t)/N_a(t)} \},
\end{eqnarray*}
where for $\mu\in\R$ and $\gamma> 0$, the set $\cC_{\mu,\gamma}$
is defined as
%
\begin{eqnarray}\label{eqCmudef}\quad
\cC_{\mu,\gamma}
&=&\bigl\{ \nu\in\mathfrak{M}_1(\cS)\dvtx  \exists\nu' \in\cD
\mbox{ with } \Ed \bigl(\nu' \bigr) > \mu\mbox{ and } \KL \bigl(
\Pi_\cD( \nu), \nu' \bigr) \leq\gamma \bigr\}.
\end{eqnarray}
By definition of $\cK_{\inf}$,
%
\begin{equation}
\label{eq6bis} \cC_{\mu,\gamma} \subseteq \bigl\{ \nu\in
\mathfrak{M}_1(\cS)\dvtx  \cK_{\inf} \bigl(\Pi_{\cD}(
\nu),\mu \bigr) \leq\gamma \bigr\}.
\end{equation}
Using (\ref{eqmajoNTa}), and recalling that for rounds $t \in\{
1,\ldots,K \}$,
each arm is played once, one obtains
\begin{eqnarray*}
\E \bigl[ N_a(T) \bigr] &\leq&1 + \sum_{t=K}^{T-1}
\P \bigl\{ \mu^{\dagger}\geq U_{a^\star}(t) \bigr\}
\\
&&{}+ \sum_{t=K}^{T-1} \P\{ \widehat{
\nu}_{a,N_a(t)} \in\cC_{\mu^{\dagger}, f(t)/N_a(t)} \mbox{ and } A_{t+1} = a \}.
\end{eqnarray*}
The two sums in this decomposition are handled separately.
The first sum is negligible with respect to the second sum:
case-specific arguments, given in Sections~\ref{secklUCB} and
\ref{secKLemp}, prove the following statement.
%
\begin{toprove}
\label{fact1}
For proper choices of $\Pi_\cD$, $f$ and $\mu^{\dagger}$,
the sum $\sum\P\{ \mu^{\dagger}\geq U_{a^\star}(t) \}$ is
negligible with
respect to $\log T$.
\end{toprove}

The second sum is thus the leading term in the bound.
It is first rewritten using the stopping times $\tau_{a,2}, \tau
_{a,3}, \ldots$
introduced in Section~\ref{secnot}. Indeed, $A_{t+1} = a$ happens for
$t \geq K$ if and only
if $\tau_{a,n} = t+1$ for some $n \in\{ 2,\ldots, t+1 \}$; and of
course, two stopping
times $\tau_{a,n}$ and $\tau_{a,n'}$ cannot be equal when $n \ne n'$.
We also note that
$N_a(\tau_{a,n}-1) = n-1$ for $n \geq2$. Therefore,
%
\begin{eqnarray}
\label{eqdecompNat}
&& \sum_{t=K}^{T-1}
\P\{ \widehat{\nu}_{a,N_a(t)} \in\cC_{\mu^{\dagger}, f(t)/N_a(t)} \mbox{ and } A_{t+1}
= a \}
\nonumber
\\
&&\qquad \leq \sum_{t=K}^{T-1} \P\{ \widehat{
\nu}_{a,N_a(t)} \in\cC_{\mu^{\dagger}, f(T)/N_a(t)} \mbox{ and } A_{t+1} = a \}
\nonumber\\
&&\qquad = \sum_{t=K}^{T-1} \sum
_{n=2}^{T-K+1} \P\{ \widehat{\nu}_{a,N_a(t)} \in
\cC_{\mu^{\dagger}, f(T)/N_a(t)} \mbox{ and } \tau_{a,n} = t+1 \}
\\
&&\qquad = \sum_{n=2}^{T-K+1} \sum
_{t=K}^{T-1} \P\{ \widehat{\nu}_{a,n-1} \in
\cC_{\mu^{\dagger}, f(T)/(n-1)} \mbox{ and } \tau_{a,n} = t+1 \}
\nonumber\\
&&\qquad \leq \sum_{n=1}^{T-K} \P\{ \widehat{
\nu}_{a,n} \in\cC_{\mu^{\dagger}, f(T)/n} \},\nonumber
\end{eqnarray}
where we used, successively, the following facts: the sets $\cC_{\mu
^{\dagger},\gamma}$ grow with~$\gamma$; the event $\{A_{t+1} = a\}$ can
be written as a disjoint union of the events $\{ \tau_{a,n} = t+1\}$,
for ${2\leq n\leq T-K+1}$; the events $\{ \tau_{a,n} = t+1 \}$ are
disjoint as $t$ varies between $K$ and $T-1$, with a possibly empty
union (as $\tau_{a,n}$ may be larger than $T$).

By upper bounding the first
%
\begin{equation}
\label{eqdefs0} n_0 = \biggl\lceil\frac{f(T)}{\cK_{\inf} (\nu
_a,\mu^\star)} \biggr\rceil
\end{equation}
terms of the sum in (\ref{eqdecompNat}) by $1$, we obtain
\[
\sum_{n=1}^{T-K} \P\{ \widehat{
\nu}_{a,n} \in\cC_{\mu^{\dagger}, f(T)/n} \} \leq\frac{f(T)}{\cK
_{\inf}
(\nu_a,\mu^\star)} + 1 + \sum
_{n \geq n_0 + 1} \P\{ \widehat{\nu}_{a,n} \in
\cC_{\mu^{\dagger}, f(T)/n} \}.
\]
It remains to upper bound the remaining sum: this is the object of the
following statement, which will also be proved using case-specific arguments.

\begin{toprove}
\label{fact3}
For proper choices of $\Pi_\cD$, $f$ and $\mu^{\dagger}$, the sum
$\sum\P\{ \widehat{\nu}_{a,n} \in\cC_{\mu^{\dagger}, f(T)/n} \}$
is negligible with respect to $\log T$.
\end{toprove}

Putting everything together, one obtains
%
\begin{eqnarray}
\label{eqmajN}\qquad \E \bigl[ N_a(T) \bigr] &\leq&\frac{f(T)}{\cK_{\inf
} (\nu_a,\mu^\star)}
\nonumber\\[-8pt]\\[-8pt]
&&{}+ \underbrace{\sum_{n \geq n_0 + 1} \P\{ \widehat{
\nu}_{a,n} \in\cC_{\mu^{\dagger}, f(T)/n} \}}_{o(\log T)}{} + {}\underbrace {\sum
_{t=K}^{T-1} \P \bigl\{ \mu^{\dagger}\geq
U_{a^\star}(t) \bigr\}}_{o(\log T)} +\, 2.
\nonumber
\end{eqnarray}
Theorems~\ref{thmklUCB} and~\ref{thmKLempUCB} are instances of this
general bound providing nonasymptotic controls for $\E[ N_a(T) ]$ in
the two settings considered in this paper.

\section{Rewards in a canonical one-dimensional exponential family}
\label{secklUCB}

We consider in this section the case when $\cD$ is a canonical
exponential family of probability distributions $\nu_\theta$,
indexed by $\theta\in\Theta$; that is, the distributions $\nu
_\theta$ are absolutely
continuous with respect to a dominating measure $\rho$ on $\R$, with
probability density
\[
\frac{\d\nu_\theta}{\d\rho}(x) = \exp \bigl( x\theta- b(\theta ) \bigr),\qquad
x \in\R;
\]
we assume in addition that $b\dvtx  \Theta\to\R$ is twice
differentiable. We also assume that
$\Theta\subseteq\R$ is the natural parameter space, that is, the set
\[
\Theta= \biggl\{ \theta\in\R\dvtx  \int_{\R} \exp(x\theta) \,\d
\rho(x) < \infty \biggr\}
\]
and that the exponential family $\cD$ is regular, that is, that
$\Theta$ is an
open interval (an assumption that turns out to be true in all the examples
listed below). In this setting, considered in the pioneering papers
by \citet{LaiRo85} and \citet{Agrawal95}, the upper
confidence bound defined
in (\ref{eqUCBdef}) takes an explicit form related to the large deviation
rate function. Indeed, as soon as the reward distributions satisfy Chernoff-type
inequalities, these can be used to construct an UCB policy, while for
heavy-tailed
distributions other approaches are required, as surveyed by \citet
{bubeckcesa-bianchi2012-ftml}.

For a thorough
introduction to canonical exponential families, as well as proofs of
the following
properties, the reader is referred to \citet{lehmanncasella98}.
The derivative $\dot{b}$ of $b$ is an increasing continuous function
such that $\Ed(\nu_\theta) = \dot{b}(\theta)$
for all $\theta\in\Theta$; in particular, $b$ is strictly convex.
Thus, $\dot{b}$ is one-to-one with a continuous inverse $\dot
{b}^{-1}$ and
the distributions $\nu_\theta$ of $\cD$ can also be parameterized by
their expectations $\Ed(\nu_\theta)$.
Defining the open interval of all expectations, $I = \dot{b}(\Theta)
= (\mu_-, \mu_+)$,
there exists a unique distribution of $\cD$ with expectation $\mu\in
I$, namely, $\nu_{\dot{b}^{-1}(\mu)}$.

The Kullback--Leibler divergence between two distributions $\nu_\theta,
\nu_{\theta'} \in\cD$ is given by
\[
\KL(\nu_\theta,\nu_{\theta'}) = \bigl(\theta-\theta'
\bigr) \dot{b}(\theta) - b(\theta) + b \bigl(\theta' \bigr),
\]
which, writing $\mu= \Ed(\nu_\theta)$ and $\mu' = \Ed(\nu
_{\theta'})$, can be reformulated as
%
\begin{eqnarray}
\label{eqKLexpfam}
d \bigl(\mu,\mu' \bigr) &\eqdef&\KL(
\nu_{\theta},\nu_{\theta'}) \nonumber\\[-8pt]\\[-8pt]
&=& \bigl( \dot{b}^{-1}(\mu)-
\dot{b}^{-1} \bigl(\mu' \bigr) \bigr) \mu- b \bigl(
\dot{b}^{-1}(\mu) \bigr) + b \bigl( \dot{b}^{-1} \bigl(
\mu' \bigr) \bigr).\nonumber
\end{eqnarray}
This defines a divergence $d\dvtx  I\times I \to\R_+$ that inherits from
the Kullback--Leibler divergence
the property that $d(\mu,\mu') = 0$ if and only if $\mu= \mu'$.
In addition, $d$ is (strictly) convex and differentiable over $I\times I$.

As the examples below of specific canonical exponential families
illustrate, the closed-form expression for this
re-parameterized Kullback--Leibler divergence is usually simple.\vadjust{\goodbreak}

\begin{example}[(Binomial distributions for $n$-samples)] $\theta= \log
(\mu/(n-\mu) )$, $\Theta= \R$, $b(\theta) = n\log(1+\exp(\theta
) )$, $I = (0,n)$,
\[
d \bigl(\mu,\mu' \bigr) = \mu\log\frac{\mu}{\mu'} + (n-\mu) \log
\frac{n-\mu}{n-\mu'}.
\]
The case $n=1$ corresponds to Bernoulli distributions.
\end{example}

\begin{example}[(Poisson distributions)]
$\theta= \log(\mu)$, $\Theta= \R$, $b(\theta)=\exp(\theta)$, $I
= (0,+\infty)$,
\[
d \bigl(\mu,\mu' \bigr) = \mu' - \mu+ \mu\log
\frac{\mu}{\mu'}.
\]
\end{example}

\begin{example}[(Negative binomial distributions with known shape
parameter $r$)] $\theta= \log(\mu/(r+\mu))$, $\Theta= (-\infty,
0)$, $b(\theta) = -r\log(1-\exp(\theta) )$,
$I = (0,+\infty)$,
\[
d \bigl(\mu,\mu' \bigr) = r \log\frac{r+\mu'}{r+\mu} + \mu\log
\frac{\mu(r+\mu')}{\mu'(r+\mu)}.
\]
The case $r=1$ corresponds to geometric distributions.
\end{example}

\begin{example}[(Gaussian distributions with known variance $\sigma^2$)]
$\theta= \mu/\sigma^2$, $\Theta= \R$, $b(\theta) = \sigma
^2\theta^2/2$, $I=\R$,
\[
d \bigl(\mu,\mu' \bigr) = \frac{(\mu- \mu')^2}{2 \sigma^2}.
\]
\end{example}

\begin{example}[(Gamma distributions with known shape parameter $\alpha$)]
$\theta= -\alpha/\mu$, $\Theta= (-\infty, 0)$, $b(\theta) =
-\alpha\log(-\theta)$, $I = (0,+\infty)$,
\[
d \bigl(\mu,\mu' \bigr) = \alpha \biggl( \frac{\mu}{\mu'} - 1 -
\log\frac{\mu}{\mu'} \biggr).
\]
The case $\alpha= 1$ corresponds to exponential distributions.
\end{example}

For all $\mu\in I$ the convex functions $d( \cdot, \mu)$ and $d(\mu,
\cdot)$ can be extended by continuity to $\oI= [\mu_-,\mu_+]$ as
follows:
\[
d(\mu_-, \mu) = \lim_{\mu'\to\mu_-} d \bigl(\mu',\mu
\bigr),\qquad d(\mu_+, \mu) = \lim_{\mu'\to\mu_+} d \bigl(\mu',
\mu \bigr)
\]
with similar statements for the second function. Note that these limits
may equal $+\infty$;
the extended function $d\dvtx  \oI\times I \cup I \times\oI\to
[0,+\infty]$
is still a convex function. By convention,
we also define $d(\mu_-, \mu_-) = d(\mu_+, \mu_+) = 0$.

Note that our exponential family models are minimal in the sense of
Wainwright and Jordan [(\citeyear{wainwrightjordan2008}), Section 3.2]
and thus that $I$ coincides with the interior of the set of realizable
expectations for all distributions that are absolutely continuous with
respect to $\rho$; see \citet{wainwrightjordan2008}, Theorem
3.3 and Appendix~B. In particular, this
implies that distributions in $\cD$ have supports in $\oI$ and that,
consequently, the
empirical means $\widehat{\nu}_a(t)$ are in $\oI$ for all $a$ and $t$.
(Note, however, that they may not be in $I$ itself: think in particular
of the case of Bernoulli
distributions when $t$ is small.)

\subsection{The \textit{\texttt{kl-UCB}} algorithm}

As the distributions in $\cD$ can be parameterized by their
expectation, $\Pi_\cD$ associates with each $\nu\in\measg(\oI)$
such that $\Ed(\nu) \in I$ the distribution $\nu_{\dot{b}^{-1}(\Ed
(\nu))} \in\cD$ which has the same expectation.

As shown above, for all $\nu' \in\cD$, it then holds that $\KL(
\Pi_\cD(\nu), \nu' )= d ( \Ed(\nu),\break \Ed(\nu') )$; and this
equality can be extended to the case where $\Ed(\nu)\in\oI$. In
this setting,
sufficient statistics for $\widehat{\nu}_a(t)$ and $\widehat{\nu}_{a,n}$ are
given by,
respectively,
\[
\widehat{\mu}_a(t) = \frac{1}{N_a(t)} \sum
_{s=1}^t Y_s \ind_{ \{ A_s = a \} }
\quad\mbox{and}\quad \widehat{\mu}_{a,n} = \frac{1}{n} \sum
_{k=1}^n X_{a,k},
\]
where the former is defined as soon as $N_a(t) \geq1$.

The upper-confidence bound $U_a(t)$ may be defined in this model
not only in terms of $\cD$ but also of its ``boundaries,''
namely, in terms of $\oI$ and not only $I$, as
%
\begin{equation}
\label{eqUCBexp} U_a(t) = \sup \biggl\{ \mu\in\oI\dvtx  d \bigl( \widehat{
\mu}_{a}(t), \mu \bigr) \leq\frac{f(t)}{N_{a}(t)} \biggr\}.
\end{equation}
This supremum is achieved: in the case when $\widehat{\mu}_{a}(t) \in I$,
this follows
from the fact that $d$ is continuous on $I \times\oI$; when $\widehat{\mu
}_{a}(t) = \mu_+$, this
is because $U_a(t) = \mu_+$; in the case when $\widehat{\mu}_{a}(t) = \mu
_-$, either $\mu_-$ is
the only $\mu\in\oI$ for which $d(\mu_-,\mu)$ is finite, or $d(\mu
_-, \cdot)$ is convex thus
continuous on the open interval where it is finite.

Thus, in the setting of this section, Algorithm~\ref{algKLUCB}
rewrites as Algorithm~\ref{algklUCB},
which will be referred to as \texttt{kl-UCB}.

\begin{algorithm}[b]
\begin{minipage}{338pt}
\caption{The \texttt{kl-UCB} algorithm.}
\label{algklUCB}

\textit{Parameters}: A nondecreasing function $f\dvtx  \mathbb{N} \to
\mathbb{R}$ \\
\textit{Initialization}: Pull each arm of $\{ 1,\ldots, K \}$ once
\\
\For{$t=K$ to $T-1$,}{

compute for each arm $a$ the quantity
\[
U_a(t) = \sup \biggl\{ \mu\in\oI\dvtx  d \bigl( \widehat{
\mu}_{a}(t), \mu \bigr) \leq\frac{f(t)}{N_{a}(t)} \biggr\}
\]
pick an arm $ {A_{t+1} \in\argmax_{a \in\{ 1,\ldots,K \}} U_a(t)}$}
\end{minipage}
\end{algorithm}

In practice, the computation of $U_a(t)$ boils down to finding the zero
of an increasing
and convex scalar function. This can be done either by dichotomic
search or by Newton iterations. In all the examples given above,
well-known inequalities
(e.g., Hoeffding's inequality) may be used to obtain an initial upper
bound on~$U_a(t)$.

\subsection{Regret analysis}
In this parametric context we have $\cK_{\inf} (\nu,\mu) = d (\Ed
(\nu),\mu)$ when $\Ed(\nu) \in I$ and $\mu\in I$.
In light of the results by \citet{LaiRo85} and \citet{Agrawal95},
the following theorem thus proves the asymptotic optimality of the
\texttt{kl-UCB} algorithm.
Moreover, it provides an explicit, nonasymptotic bound on the regret.

\begin{theorem}
\label{thmklUCB}
Assume that all arms belong to a canonical, regular, exponential family
$\mathcal{D} = \{\nu_\theta\dvtx  \theta\in\Theta\}$ of probability
distributions indexed by its natural parameter space $\Theta\subseteq
\R$. Then, using Algorithm~\ref{algklUCB} with the divergence $d$
given in (\ref{eqKLexpfam}) and with the choice $f(t) = \log
(t)+3\log\log(t)$ for $t \geq3$ and $f(1) = f(2) = f(3)$, the number
of draws of any suboptimal arm $a$ is upper bounded for any horizon $T
\geq3$ as
\begin{eqnarray*}
\E \bigl[N_a(T) \bigr] &\leq&\frac{\log(T)}{d (\mu_a,\mu^\star)} + 2\sqrt{
\frac{2\pi\sigma_{a,\star}^2 ({d'}(\mu_a, \mu^\star) )^2}{ (
d(\mu_a,\mu^\star) )^3}} {\sqrt{\log(T) + 3\log \bigl(\log (T) \bigr)}}
\\
&&{}+ \biggl(4e + \frac{3}{d(\mu_a,\mu^\star)} \biggr)\log \bigl(\log(T) \bigr) + 8
\sigma_{a,\star}^2 \biggl( \frac{d'(\mu_a,\mu^\star)}{d(\mu
_a,\mu^\star)}
\biggr)^2 + 6,
\end{eqnarray*}
where $\sigma_{a,\star}^2 = \max\{\Vard(\nu_\theta)\dvtx  \mu_a
\leq\Ed(\nu_\theta) \leq\mu^\star\}$ and where $d'( \cdot, \mu
^\star)$ denotes the derivative of $d( \cdot, \mu^\star)$.
\end{theorem}

The proof of this theorem is provided in the supplemental
article [\citet{extended}, Appendix A]. A key argument, proved in
Lemma 2 (see also Lemma 11), is the following deviation bound for the
empirical mean with random number of summands: for all $\epsilon>1$
and all $t\geq1$,
%
\begin{equation}
\label{eqinegDevIndepInterest}\quad 
\P \biggl\{
\widehat{\mu}_{a^\star}(t) < \mu^\star \mbox{ and } d \bigl(\widehat{
\mu}_{a^\star}(t), \mu^\star \bigr) \geq\frac{\epsilon}{N_{a^\star}(t)} \biggr\}
\leq e \bigl\lceil \epsilon\log(t) \bigr\rceil\exp(-\epsilon).
\end{equation}

For binary distributions, guarantees analogous to that of Theorem \ref
{thmklUCB} have been obtained recently for algorithms inspired by the Bayesian
paradigm, including the so-called \citet{Tho33} sampling
strategy, which is not an index policy in the sense of \citet
{Agrawal95}; see
\citet{12-kaufmann-cappe-garivieraistats} and \citet{EKALT-12}.

\section{Bounded and finitely supported rewards}
\label{secKLemp}

In this section, $\cD$ is the set $\cF$ of finitely supported
probability distributions
over $\cS= [0,1]$. In this case, the empirical measures $\widehat{\nu
}_a(t)$ belong to $\cF$ and
hence the operator $\Pi_\cD$ is taken to be the identity. We denote
by $\Supp(\nu)$ the
finite support of an element $\nu\in\cF$.\vadjust{\goodbreak}

The maximization program (\ref{eqUCBdef}) defining $U_a(t)$ admits
in this case
the simpler formulation
\begin{eqnarray*}
U_a(t) &\eqdef&\sup \biggl\{ \Ed(\nu)\dvtx  \nu\in\cF\mbox{ and } \KL
\bigl( \widehat{\nu}_a(t), \nu \bigr) \leq\frac{f(t)}{N_a(t)} \biggr\}
\\
&=& \sup \biggl\{ \Ed(\nu)\dvtx  \nu\in\mathfrak{M}_1 \bigl(\Supp \bigl(
\widehat{\nu}_a(t) \bigr) \cup\{1\} \bigr) \mbox{ and } \KL \bigl( \widehat{
\nu}_a(t), \nu \bigr) \leq\frac{f(t)}{N_a(t)} \biggr\},
\end{eqnarray*}
which admits an explicit computational solution;
these two points are detailed in the supplemental article [\citet
{extended}, Appendix C.1].
The reasons for which the value 1 needs to be added to the support (if
it is not yet present) will be detailed in Section~\ref{secKLUCBgeneral}.

Thus Algorithm~\ref{algKLUCB} takes the following simpler form, which
will be referred to as
the empirical \texttt{KL-UCB} algorithm.

Like the \texttt{DMED} algorithm, for which asymptotic bounds are
proved in Honda and Takemura (\citeyear{HondaTakemura10DMED,HondaTakemura11ML}),
Algorithm~\ref{algKLUCB} relies on the empirical likelihood method
[see \citet{owen2001empirical}] for the construction of the
confidence bounds. However, \texttt{DMED} is not an index policy, but
it maintains a list of active arms---an approach that, generally
speaking, seems to be less satisfactory and slightly less efficient in
practice. Besides, the analyses of the two algorithms, even though they
both rely on some technical properties of the function $\cK_{\inf}$,
differ significantly.

\begin{algorithm}[b]
\begin{minipage}{340pt}
\caption{The empirical \texttt{KL-UCB} algorithm.}
\label{algKLempUCB}

\textit{Parameters}: A nondecreasing function $f\dvtx  \mathbb{N} \to
(0,+\infty)$\\
\textit{Initialization}: Pull each arm of $\{ 1,\ldots, K \}$ once
\\
\For{$t=K$ to $T-1$,}{

compute for each arm $a$ the quantity
\[
U_a(t) \!=\! \sup \biggl\{ \Ed(\nu)\dvtx  \nu\!\in\!\mathfrak{M}_1
\bigl(\Supp \bigl( \widehat{\nu}_a(t) \bigr)\! \cup\!\{1\} \bigr)
\mbox{ and }
\KL \bigl( \widehat{\nu}_a(t), \nu \bigr)\! \leq\!\frac{f(t)}{N_a(t)}\! \biggr\}
\]
pick an arm $ {A_{t+1} \in\argmax_{a \in\{ 1,\ldots,K \}} U_a(t)}$
}
\end{minipage}
\end{algorithm}

\begin{theorem}
\label{thmKLempUCB}
Assume that $\mu_a > 0$ for all arms $a$ and that $\mu^\star< 1$.
There exists a constant $M(\nu_a,\mu^\star) > 0$ only depending on
$\nu_a$ and $\mu^\star$
such that, with the choice $f(t) = \log(t) + \log(\log(t) )$ for $t
\geq2$,
the expected number of times that any suboptimal arm $a$ is pulled by
Algorithm~\ref{algKLempUCB} is smaller,
for all $T \geq3$, than
\begin{eqnarray*}
\E \bigl[ N_a(T) \bigr] &\leq& \frac{\log(T)}{\cK_{\inf} (\nu
_a,\mu^\star)} +
\frac{36}{(\mu^\star)^4} \bigl( \log(T) \bigr)^{4/5} \log \bigl(\log(T) \bigr)
\\
& &{} + \biggl( \frac{72}{(\mu^\star)^4} + \frac{2\mu^\star
}{(1-\mu^\star) \cK_{\inf} (\nu_a,\mu^\star)^2} \biggr) \bigl( \log(T)
\bigr)^{4/5}
\\
& &{} + \frac{(1-\mu^\star)^2 M(\nu_a,\mu^\star)}{2(\mu^\star)^2} \bigl( \log(T) \bigr)^{2/5}
\\
& &{} + \frac{\log(\log(T) )}{\cK_{\inf} (\nu_a,\mu^\star)} + \frac{2\mu^\star}{(1-\mu^\star) \cK_{\inf} (\nu_a,\mu^\star
)^2} + 4.
\end{eqnarray*}
\end{theorem}

Theorem~\ref{thmKLempUCB} implies a nonasymptotic bound of the form
\[
\E \bigl[ N_a(T) \bigr] \leq\frac{\log(T)}{\cK_{\inf} (\nu
_a,\mu^\star)} + O \bigl( \bigl(
\log(T) \bigr)^{4/5} \log \bigl(\log(T) \bigr) \bigr).
\]
The exact value of the constant $M(\nu_a,\mu^\star)$ is provided in
the proof of Theorem~\ref{thmKLempUCB}, which can be
found in the supplemental article [\citet{extended}, Appendix B];
see, in particular, Section B.3 as well as the variational form of $\cK
_{\inf}$ introduced in Lemma 4 of Section B.1 of the supplement.

\section{Algorithms for general bounded rewards}
\label{secgeneral}

In this section, we consider the case where the arms are only known to have
bounded distributions. As in Section~\ref{secKLemp}, we assume
without loss of
generality that the rewards are bounded in $[0,1]$. This is the setting
considered by \citet{Auer02}, where the \texttt{UCB} algorithm was
described and
analyzed. We first prove that \texttt{kl-UCB}
(Algorithm~\ref{algklUCB}) with Kullback--Leibler divergence for
Bernoulli distributions
is always preferable to \texttt{UCB}, in the sense that a smaller
finite-time regret
bound is guaranteed. \texttt{UCB} is indeed nothing but
\texttt{kl-UCB} with quadratic
divergence and we obtain a refined analysis of \texttt{UCB} as a consequence of
Theorem~\ref{thmklUCB}. We then discuss the use of the empirical
\texttt{KL-UCB}
approach, in which one directly applies Algorithm~\ref{algKLempUCB}.
We provide preliminary
results to support the observation that empirical \texttt{KL-UCB} achieves improved
performance on sufficiently long horizons (see simulation results in
Section~\ref{secexpnum}), at the price, however, of a
significantly higher computational complexity.

\subsection{The \textit{\texttt{kl-UCB}} algorithm for bounded
distributions}
A careful reading of the proof of Theorem~\ref{thmklUCB} (see the
supplemental article [\citet{extended}, Section~A])
shows that \texttt{kl-UCB} enjoys regret guarantees
in models with arbitrary bounded distributions $\nu$ over $[0,1]$ as
long as it is
used with a divergence $d$ over $[0,1]^2$ satisfying the following
double property:
there exists a family of strictly convex and continuously
differentiable functions $\phi_\mu\dvtx  \R\to[0,+\infty)$,
indexed by $\mu\in[0,1]$, such that
first, $d( \cdot,\mu)$ is the convex conjugate of $\phi_{\mu}$ for
all $\mu\in[0,1]$; and, second,
the domination condition $\mathcal{L}_{\nu}(\lambda) \leq\phi_{\Ed
(\nu)}(\lambda)$ for all $\lambda\in\R$ and all $\nu\in\measg(
[0,1] )$ holds, where $\mathcal{L}_{\nu}$ denotes the\vadjust{\goodbreak}
moment-generating function of $\nu$,
\[
\mathcal{L}_{\nu}\dvtx  \lambda\in\R\longmapsto\mathcal{L}_{\nu}(
\lambda) = \int_{[0,1]} e^{\lambda x} \,\d\nu(x).
\]

The following elementary lemma dates back to \citet{Hoeffding63};
it upper bounds the moment-generating function of any probability
distribution over $[0,1]$
with expectation $\mu$ by the moment-generating function of the
Bernoulli distribution with parameter $\mu$,
which is further bounded by the moment-generating function of the
normal distribution with mean $\mu$ and variance $1/4$.
All these moment-generating functions are defined on the whole real
line $\R$.
In light of the above, it thus shows that the
Kullback--Leibler divergence $d_{\ber}$ between Bernoulli
distributions and the
Kullback--Leibler divergence $d_{\qua}$ between normal distributions
with variance $1/4$
are adequate candidates for use in the \texttt{kl-UCB}
algorithm in the case of bounded distributions.
%
\begin{lemma}
\label{lembounded2bernoulli}
Let $\nu\in\measg( [0,1] )$ and let $\mu= \Ed(\nu)$. Then, for
all $\lambda\in\R$,
\[
\mathcal{L}_{\nu}(\lambda) = \int_{[0,1]}
e^{\lambda x} \,\d\nu(x) \leq1-\mu+ \mu\exp(\lambda) \leq\exp \bigl( \lambda
\mu+ 2\lambda^2 \bigr).
\]
\end{lemma}

The proof of this lemma is straightforward; the first inequality is by
convexity,
as $e^{\lambda x} \leq x e^{\lambda} + (1-x)$ for all $x \in[0,1]$,
and the second inequality
follows by standard analysis.

We therefore have the following corollaries to Theorem
\ref{thmklUCB}. (They are obtained by
bounding in particular the variance term $\sigma^2_{a,\star}$ by $1/4$.)

\begin{corollary}\label{thborneRegretklUCB}
Consider a bandit problem with rewards bounded in $[0,1]$.
Choosing the parameters $f(t) = \log(t) + 3\log\log(t)$ for $t \geq
3$ and $f(1) = f(2) = f(3)$, and
\[
d_{\ber} \bigl(\mu,\mu' \bigr) = \mu\log\frac{\mu}{\mu'}
+ (1-\mu)\log\frac{1-\mu}{1-\mu'}
\]
in Algorithm~\ref{algklUCB}, the number of draws of any suboptimal
arm $a$ is
upper bounded for any horizon $T \geq3$ as
\begin{eqnarray*}
\E \bigl[N_a(T) \bigr] &\leq&\frac{\log(T)}{d_{\ber}(\mu_a,\mu
^\star)} \\
&&{}+
\frac{\sqrt{2\pi} \log({\mu^ \star(1-\mu_a)}/({\mu_a(1-\mu
^\star)}) )}{
( d_{\ber}(\mu_a,\mu^\star) )^{3/2}} {\sqrt{\log(T) + 3\log \bigl(\log(T) \bigr)}}
\\
&&{}+ \biggl(4e + \frac{3}{d_{\ber}(\mu_a,\mu^\star)} \biggr) \log \bigl(\log(T) \bigr)\\
&&{} +
\frac{2 ( \log({\mu^ \star(1-\mu_a)}/({\mu_a(1-\mu^\star)})
) )^2}{ (d_{\ber}(\mu_a,\mu^\star) )^2} + 6.
\end{eqnarray*}
\end{corollary}

We denote by $\phi_{\Ed(\nu)} = 1 - \Ed(\nu) + \Ed(\nu)\exp(
\cdot)$ the upper bound
on $\mathcal{L}_{\nu}$ exhibited in Lemma
\ref{lembounded2bernoulli}. Standard
results on Kullback--Leibler divergences are
that for all $\mu,\mu' \in[0,1]$ and all $\nu,\nu' \in\measg(
[0,1] )$,
\[
d_{\ber} \bigl(\mu,\mu' \bigr) = \sup_{\lambda\in\R}
\bigl\{ \lambda\mu- \phi_{\mu'}(\lambda) \bigr\}
\quad\mbox{and}\quad \KL \bigl(
\nu, \nu' \bigr) \geq\sup_{\lambda\in\R} \bigl\{ \lambda
\Ed(\nu)-\mathcal{L}_{\nu'}(\lambda) \bigr\};
\]
see Massart [(\citeyear{massart2007}), pages 21 and 28]; see also \citet{DeZe98}.
Because of Lemma~\ref{lembounded2bernoulli}, it thus holds that for
all distributions $\nu,\nu' \in\measg( [0,1] )$,
\[
d_{\ber} \bigl( \Ed(\nu), \Ed \bigl(\nu' \bigr) \bigr) \leq
\KL \bigl(\nu,\nu' \bigr),
\]
and it follows that in the model $\cD= \measg( [0,1] )$ one has
\[
\cK_{\inf} \bigl(\nu_a,\mu^\star \bigr) \geq
d_{\ber} \bigl(\mu_a,\mu^\star \bigr).
\]
As expected, the \texttt{kl-UCB} algorithm may not be
optimal for all sub-families of bounded distributions.
Yet, this algorithm has stronger guarantees than the \texttt{UCB} algorithm.
It is readily checked that the latter exactly corresponds to the choice of
\[
d_{\qua} \bigl(\mu,\mu' \bigr) = 2 \bigl(\mu-
\mu' \bigr)^2
\]
in Algorithm~\ref{algklUCB} together with some nondecreasing
function $f$. For instance,
the original algorithm \texttt{UCB1} of Auer, Cesa-Bianchi and Fischer
[(\citeyear{Auer02}), Theorem 1], relies on $f(t) = 4 \log(t)$.
The analysis derived in this paper gives an improved analysis of the
performance of the \texttt{UCB} algorithm
by resorting to the function $f$ described in the statement of
Theorem~\ref{thmklUCB}.
%
\begin{corollary}
\label{corUCB}
Consider the \textit{\texttt{kl-UCB}} algorithm with $d_{\qua}$
and the function $f$ defined in Theorem~\ref{thmklUCB},
or equivalently, the \textit{\texttt{UCB}} algorithm tuned as follows: at step
$t+1>K$, an arm maximizing the
upper-confidence bounds
\[
\widehat{\mu}_a(t) + \sqrt{ \bigl( \log(t)+3\log\log(t) \bigr)/ \bigl(2
N_a(t) \bigr)}
\]
is chosen.
Then the number of draws of a suboptimal arm $a$ is upper bounded as
\begin{eqnarray*}
\E \bigl[N_a(T) \bigr] &\leq&\frac{\log(T)}{2(\mu^\star- \mu
_a)^2} +\frac{ 2 \sqrt{\pi}}{ (\mu^\star-\mu_a )^2} {
\sqrt{\log(T) + 3\log \bigl(\log(T) \bigr)}}
\\
&&{}+ \biggl(4e + \frac{3}{2(\mu^\star-\mu_a)^2} \biggr)\log \bigl(\log(T) \bigr) +
\frac{8}{(\mu^\star-\mu_a)^2} + 6.
\end{eqnarray*}
\end{corollary}

As claimed, it can be checked that the leading term in the bound of
Corollary~\ref{thborneRegretklUCB} is smaller than the one of Corollary
\ref{corUCB} by applying Pinsker's inequality $d_{\ber} \geq d_{\qua}$.
The bound obtained in Corollary~\ref{corUCB} above also improves on the
one of Auer, Cesa-Bianchi and Fischer [(\citeyear{Auer02}), Theorem~1],
and it is ``optimal'' in the sense that the constant $1/2$ in the
logarithmic term cannot be improved. Note that a constant in front on
the leading term of the regret bound is proven to be arbitrarily close
to (but strictly greater than) $1/2$ for the \texttt{UCB2} algorithm of
\citet{Auer02}, when the parameter $\alpha$ goes to $0$ as the
horizon grows, but then other terms are unbounded. In comparison,
Corollary~\ref{corUCB} provides a bound for \texttt{UCB} with a leading
optimal constant $1/2$, and all the remaining terms of the bound are
finite and made explicit. Note, in addition, that the choice of the
parameter~$\alpha$, which drives the length of the phases during which
a single arm is played, is important but difficult in practice, where
\texttt{UCB2} does not really prove more efficient than \texttt{UCB}.

\subsection{The empirical \textit{\texttt{KL-UCB}} algorithm for bounded distributions}
\label{secKLUCBgeneral}

The justification of the use of empirical \texttt{KL-UCB} for general
bounded distributions $\mathfrak{M}_1([0,1])$ relies on the following
result.

\subsubsection*{A result of independent interest, connected to the
empirical-likelihood\break method}
The empirical-likelihood (or EL in short) method
provides a way to construct confidence bounds for the true expectation of
i.i.d. observations; for a thorough introduction to this theory,
see \citet{owen2001empirical}. We only recall briefly its
principle. Given a
sample $X_1,\ldots,X_n$ of an unknown distribution $\nu_0$, and denoting
$\widehat{\nu}_n = n^{-1}\sum_{k=1}^n \delta_{X_k}$ the empirical
distribution of
this sample, an EL upper-confidence bound for the expectation $\Ed(\nu
_0)$ of
$\nu_0$ is given by
%
\begin{equation}
\label{eqsEL} U_{\mathrm{EL}}(\widehat{\nu}_n,\epsilon) = \sup \bigl\{
\Ed \bigl(\nu' \bigr)\dvtx  \nu' \in\measg \bigl( \Supp(
\widehat{ \nu}_n) \bigr) \mbox{ and } \KL \bigl(\widehat{\nu}_n,
\nu' \bigr) \leq\epsilon \bigr\},
\end{equation}
where $\epsilon> 0$ is a parameter controlling the confidence level.

An apparent impediment to the application of this method in bandit
problems is the impossibility of obtaining nonasymptotic guarantees
for the covering probability of EL upper-confidence bounds.
In fact, it appears in (\ref{eqsEL}) that $U_{\mathrm{EL}}(\widehat{\nu
}_n,\epsilon)$ necessarily belongs to the convex envelop of the
observations. If, for example, all the observations are equal to $0$,
then $U_{\mathrm{EL}}(\widehat{\nu}_n,\epsilon)$ is also equal to $0$,
no matter what the value of $ \epsilon$ is; therefore,
it is not possible to obtain upper-confidence bounds for all confidence levels.

In the case of (upper-)bounded variables,
this problem can be circumvented by adding to the support of $\widehat{\nu
}_n$ the maximal possible value. In our case, instead of considering
$U_{\mathrm{EL}}(\widehat{\nu}_n,\epsilon)$, one should use
%
\begin{equation}
\label{eqtweakedEL}\qquad U(\widehat{\nu}_n, \epsilon) = \sup \bigl\{ \Ed
\bigl(\nu' \bigr)\dvtx  \nu' \in\measg \bigl( \Supp(\widehat{
\nu}_n) \cup\{1\} \bigr) \mbox{ and } \KL \bigl(\widehat{
\nu}_n, \nu' \bigr) \leq\epsilon \bigr\}.
\end{equation}
This idea was introduced in Honda and Takemura
(\citeyear{HondaTakemura10DMED,HondaTakemura11ML}), independently of
the EL literature. The following guarantee can be obtained; its proof
is provided in the supplemental article [\citet{extended}, Section~C.2].

\begin{proposition}
\label{propcovEL}$\!\!$
Let $\nu_0\in\mathfrak{M}_1([0,1])$ with $\Ed(\nu_0) \in(0,1)$
and let
$X_1,\ldots,X_n$ be independent random variables with common
distribution $\nu_0 \in\measg( [0,1] )$,
not necessarily with finite support.
Then, for all $\epsilon>0$,
\[
\P \bigl\{ U(\widehat{\nu}_n, \epsilon) \leq\Ed(\nu_0)
\bigr\} \leq\P \bigl\{ \cK_{\inf} \bigl(\widehat{\nu}_n, \Ed(
\nu_0) \bigr) \geq\epsilon \bigr\} \leq e(n+2)\exp(-n \epsilon),
\]
where $\cK_{\inf}$ is defined in terms of the model $\cD= \cF$.
\end{proposition}

For $\{0,1\}$-valued observations, it is readily seen that $U(\widehat
{\nu}_n, \epsilon)$ boils down to the upper-confidence bound given
by (\ref{eqUCBexp}). This example and some numerical simulations suggest
that the above proposition is not (always) optimal: the presence of the
factor $n$ in front of the exponential $\exp(-n \epsilon)$
term is indeed questionable.

\subsubsection*{Conjectured regret guarantees of empirical
\textit{\texttt{KL-UCB}}}

The analysis of empirical \texttt{KL-UCB} in the case where the arms are associated
with general bounded distributions is a work in progress. In view of
Proposition~\ref{propcovEL} and of the discussion above, it is only
the proof
of Fact \hyperref[fact3]{2} that needs to be extended.

As a preliminary result, we can prove an asymptotic regret bound, which
is indeed optimal, but for a variant of Algorithm~\ref{algKLempUCB};
it consists of playing in regimes $r$ of increasing lengths instances
of the empirical \texttt{KL-UCB} algorithm in which
the upper confidence bounds are given by
\[
\sup \biggl\{ \Ed(\nu)\dvtx \nu\in\mathfrak{M}_1 \bigl(\Supp \bigl(
\widehat{\nu}_a(t) \bigr) \cup\{1+\delta_r\} \bigr) \mbox{
and } \KL \bigl( \widehat{\nu}_a(t), \nu \bigr) \leq\frac{f(t)}{N_a(t)}
\biggr\},
\]
where $\delta_r \to0$ as the index of the regime $r$ increases.

The open questions would be to get an optimal bound for Algorithm \ref
{algKLempUCB} itself,
preferably a nonasymptotic one like those of Theorems~\ref{thmklUCB}
and~\ref{thmKLempUCB}.
Also, a computational issue arises: as the support of each empirical
distribution may contain as many points as the number of times the
corresponding arm was pulled,
the computational complexity of the empirical \texttt{KL-UCB} algorithm
grows, approximately linearly, with the number of rounds. Hence the
empirical \texttt{KL-UCB} algorithm as it stands is only suitable for small
to medium horizons (typically less than ten thousands rounds).
To reduce the numerical complexity of this algorithm without renouncing
to performance, a possible direction could be to cluster the rewards
on adaptive grids that are to be refined over time.

\section{Numerical experiments}
\label{secexpnum} The results of the previous sections show that the
\texttt{kl-UCB} and the empirical \texttt{KL-UCB} algorithms are
efficient not only in the special frameworks for which they were
developed, but also for general bounded distributions. In the rest of
this section, we support this claim by numerical experiments that
compare these methods with competitors such as \texttt{UCB} and
\texttt{UCB-Tuned}\vadjust{\goodbreak} [\citet{Auer02}], \texttt{MOSS}
[\citet{audibertbubeck2010minimax}], \texttt{UCB-V}
[\citet{AudibertEtAlUCBV}] or \texttt{DMED} [Honda and Takemura
(\citeyear{HondaTakemura10DMED,HondaTakemura11ML})]. In these
simulations, similar confidence levels are chosen for all the upper
confidence bounds, corresponding to $f(t) = \log(t)$---a choice which
we recommend in practice. Indeed, using $f(t) = \log(t) + 3\log\log(t)$
or $f(t)=(1+\epsilon)\log(t)$ (with a small $\epsilon>0$) yields
similar conclusions regarding the ranking of the performance of the
algorithms, but leads to slightly higher average regrets. More
precisely, the upper-confidence bounds we used were $U_a(t) =
\widehat{\mu}_a(t)+\sqrt{\log(t)/(2N_a(t))}$ for \texttt{UCB},
\[
U_a(t) = \widehat{\mu}_a(t)+\sqrt{\frac{2 \widehat{v}_a(t)\log(t)}{N_a(t)}}+3
\frac{\log(t)}{N_a(t)}
\]
with
%
\begin{equation}
\label{eqUCBV} \widehat{v}_a(t) = \Biggl(\frac{1}{N_a(t)} \sum
_{s=1}^t Y_s^2
\ind_{ \{ A_s = a \} } \Biggr) - \widehat{\mu}_a(t)^2
\end{equation}
for \texttt{UCB-V} and, following \citet{Auer02},
\[
U_a(t) = \widehat{\mu}_a(t)+\sqrt{\frac{\min\{1/4, \widehat
{v}_a(t)+\sqrt{2\log(t)/N_a(t)} \}\log(t)}{N_a(t)}}
\]
for \texttt{UCB-Tuned}.
Both \texttt{UCB-V} and \texttt{UCB-Tuned} are expected to improve
over \texttt{UCB} by estimating the variance of the rewards; but
\texttt{UCB-Tuned} was introduced as an heuristic improvement over
\texttt{UCB} (and does not come with a performance bound) while
\texttt{UCB-V} was analyzed by \citet{AudibertEtAlUCBV}.

Different choices of the divergence function $d$ lead to different
variants of the \texttt{kl-UCB} algorithm, which are
sometimes compared with one another in the sequel. In order to clarify
this point, we reserve the term \texttt{kl-UCB} for the
variant using the \emph{binary} Kullback--Leibler divergence
(i.e., between Bernoulli distributions), while other choices are
explicitly specified by their denomination (e.g., \texttt{kl-poisson-UCB} or
\texttt{kl-exp-UCB}
for families of Poisson or exponential distributions).
The simulations presented in this section have been performed using the
\texttt{py/maBandits} package [\citet{CapGarKau12}], which is
publicly available from the \texttt{mloss.org} website and can be used
to replicate these experiments.

\subsection{Bernoulli rewards}
We first consider the case of Bernoulli rewards, which has a special
historical importance
and which covers several important practical applications of bandit
algorithms; see \citet{ro52,gittins-1979} and references therein.
With $\{0,1\}$-valued rewards and with the binary Kullback--Leibler
divergence as a divergence function, it is readily checked that the
\texttt{kl-UCB} algorithm coincides exactly with empirical
\texttt{KL-UCB}.\vadjust{\goodbreak}

\begin{figure}

\includegraphics{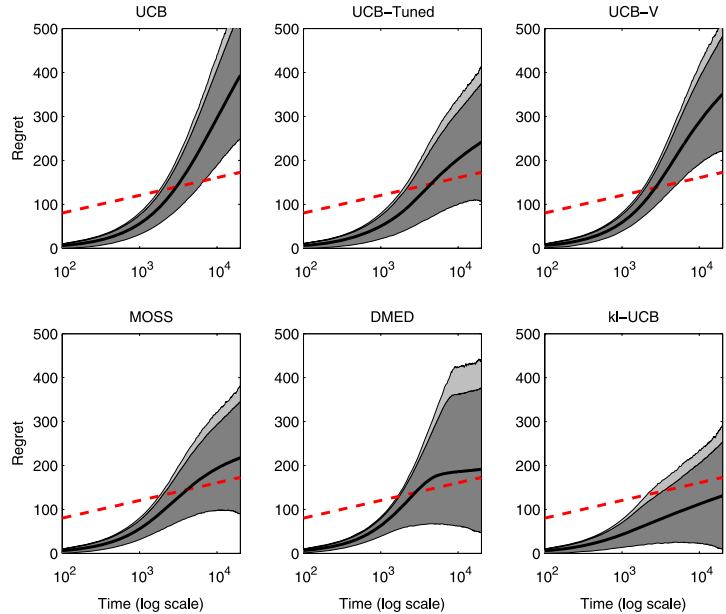}

\caption{Regret of the various algorithms as a function of time (on a
log-scale) in the
Bernoulli ten-arm scenario. On each figure, the dashed line shows the
asymptotic lower bound; the solid
bold curve corresponds to the mean regret; while the dark and light
shaded regions show,
respectively, the central $99 \%$ region and the upper $99.95 \%$ quantile.}
\label{fig10}
\end{figure}

In Figure~\ref{fig10} we consider a difficult scenario,
inspired by a situation (frequent in applications like marketing or
Internet advertising)
where the mean reward of each arm is very low. In our scenario, there
are ten arms: the
optimal arm has expected reward $0.1$, and the nine suboptimal arms
consist of three
different groups of three (stochastically) identical arms, each with
respective expected rewards $0.05$,
$0.02$ and $0.01$. We resorted to $N=50\mbox{,}000$ simulations to obtain the regret
plots of Figure~\ref{fig10}. These plots show, for each algorithm,
the average cumulated
regret together with quantiles of the cumulated regret distribution as
a function of time
(on a logarithmic scale).

Here, there is a huge gap in performance between \texttt{UCB} and
\texttt{kl-UCB}. This is explained by the fact
that the variances of all reward
distributions are much smaller than $1/4$, the pessimistic upper bound
used in
Hoeffding's inequality (i.e., in the design of \texttt{UCB}). The gain in
performance of
\texttt{UCB-Tuned} is not very significant.
\texttt{kl-UCB} and \texttt{DMED}
reach a performance that is on par with the lower bound (\ref{eqbinfRegret})
of \citet{Burnetas96} (shown in strong dashed line); the
performance of
\texttt{kl-UCB} is somewhat better than the one of
\texttt{DMED}. Notice that for the best methods, and in
particular for \texttt{kl-UCB}, the mean regret is below
the lower bound,
even for larger horizons, which reveals and illustrates the asymptotic
nature of this bound.

\subsection{Truncated Poisson rewards}

In this second scenario, we consider $6$ arms with truncated Poisson
distributions. More precisely, each arm $1\leq a\leq6$ is associated
with $\nu_a$, a Poisson
distribution with expectation $(2+a)/4$, truncated at $10$. The
experiment consisted of
$N=10\mbox{,}000$ Monte Carlo replications on an horizon of $T=20\mbox{,}000$ steps. Note
that the truncation does not alter much the distributions here, as the
probability of draws larger than $10$ is small for all arms. In fact, the
role of this truncation is only to provide an explicit upper bound on the
possible rewards, which is required for most algorithms.

\begin{figure}

\includegraphics{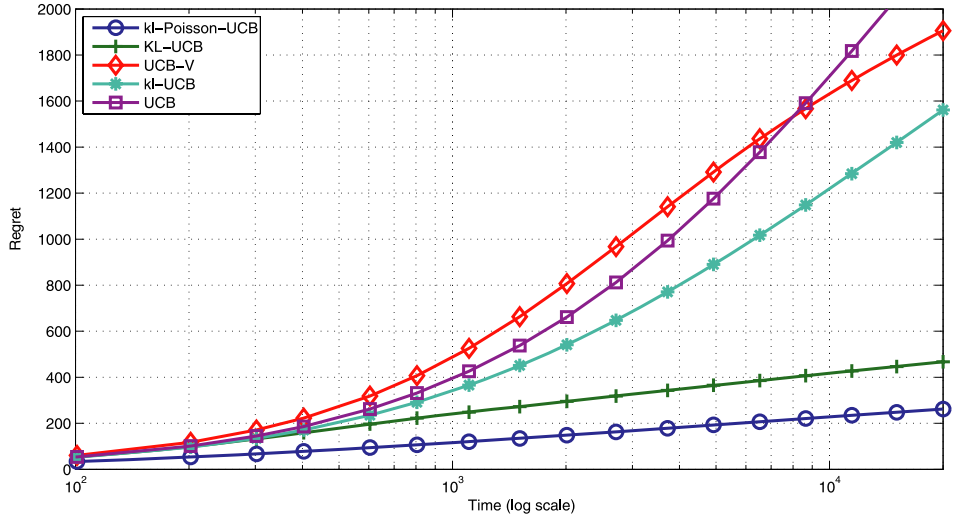}

\caption{Regret of the various algorithms as a function of time in the
truncated Poisson
scenario.}
\label{figtruncPoisson}
\end{figure}

Figure~\ref{figtruncPoisson} shows that, in this case again, the
\texttt{UCB}
algorithm is significantly worse than some of its competitors. The
\texttt{UCB-V} algorithm, which appears to have a larger regret on the first
$5000$ steps, progressively improves thanks to its use of variance
estimates for
the arms. But the horizon $T=20\mbox{,}000$ is (by far) not sufficient for
\texttt{UCB-V} to provide an
advantage over \texttt{kl-UCB}, which is thus seen to
offer an interesting
alternative even in nonbinary cases.

These three methods, however, are outperformed by the \texttt{kl-poisson-UCB}
algorithm: using the properties of the Poisson distributions (but not
taking truncation into account, however), this algorithm achieves a regret
that is about ten times smaller. In-between stands the empirical \texttt{KL-UCB}
algorithm;
it relies on nonparametric empirical-likelihood-based upper bounds and
is therefore
distribution-free as explained in Section~\ref{secKLUCBgeneral}, yet,
it proves remarkably efficient.

\subsection{Truncated exponential rewards}
In the third and last example, there are $5$ arms associated with
continuous distributions: the
rewards are exponential variables, with respective parameters $1/5$,
$1/4$, $1/3$, $1/2$
and $1$, truncated at $x_{\max}=10$ (i.e., they are bounded in $[0,10]$).

\begin{figure}

\includegraphics{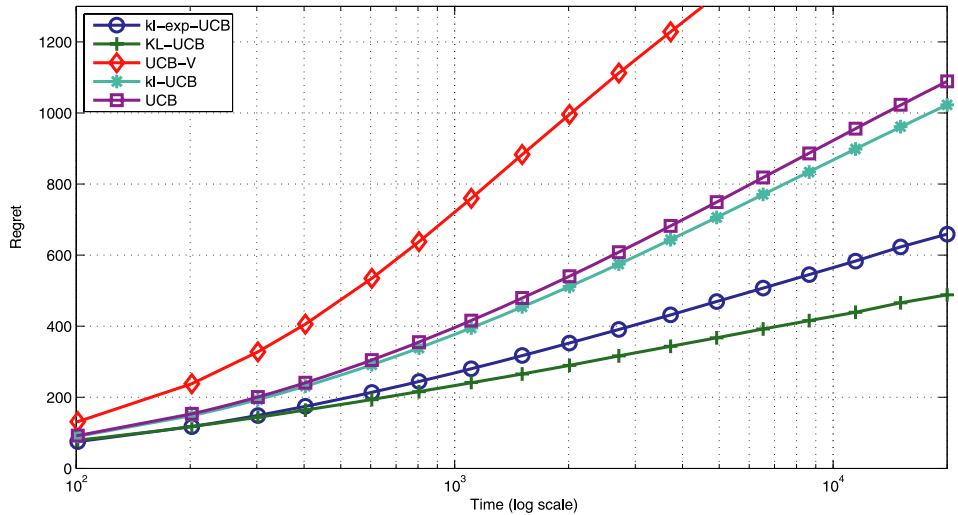}

\caption{Regret of the various algorithms as a function of time in the
truncated
exponential scenario.}
\label{figEB5}
\end{figure}

Figure~\ref{figEB5} shows that in this scenario,
\texttt{UCB} and \texttt{MOSS} are clearly
suboptimal. This time, the
\texttt{kl-UCB} does not provide a significant improvement
over \texttt{UCB} as the expectations
of the arms are not particularly close to $0$ or to $x_{\max}=10$;
hence the
confidence intervals computed by \texttt{kl-UCB} are close
to those used by \texttt{UCB}.
\mbox{\texttt{UCB-V}}, by estimating the variances of the distributions of
the rewards, which
are much smaller than the variances of $\{0,10\}$-valued distributions
with the
same expectations, would be expected to perform significantly better.
But here
again, \texttt{UCB-V} is not competitive, at least for a horizon
$T=20\mbox{,}000$. This can be explained by the fact that the upper confidence
bound of any suboptimal arm $a$, as stated in~(\ref{eqUCBV}),
contains a residual term $3\log(t) / N_a(t)$; this term is negligible
in common applications of Bernstein's inequality, but it does not
vanish here because $N_a(t)$ is precisely of order $\log(t)$; see also
\citet{GaCa11} for further discussion of this issue.

The \texttt{kl-exp-UCB}
algorithm uses the divergence $d(x,y)=x/y-1-\log(x/y)$ prescribed for
genuine exponential distributions, but it ignores the fact that the
rewards are
truncated. However, contrary to the previous scenario, the truncation
has an
important effect here, as values larger than $10$ are relatively
probable for each arm. Because \texttt{kl-exp-UCB} is not aware of
the truncation, it uses upper bounds that are slightly too large;
however, the
performance is still excellent and stable, and the algorithm is particularly
simple.

But the best-performing algorithm in this case is the nonparametric
algorithm, empirical \texttt{KL-UCB}. This method
appears to reach here the best compromise between efficiency and
versatility, at the
price of a larger computational complexity.

\section{Conclusion}
The \texttt{kl-UCB} algorithm is a quasi-optimal method
for multi-armed bandits whenever the distributions associated with the
arms are known to belong to a simple parametric family. For each
one-dimensional exponential family, a specific divergence function has
to be used in order to achieve the lower bound (\ref{eqbinfRegret})
of \citet{LaiRo85}.

However, the binary Kullback--Leibler divergence plays a special role:
it is a conservative, universal choice for bounded distributions. The
resulting algorithm is versatile, fast and simple and proves to be a
significant improvement, both in theory and in practice, over the
widely used \texttt{UCB} algorithm.

The more elaborate \texttt{KL-UCB} algorithm relies on nonparametric
inference, by using the so-called empirical likelihood method. It is
optimal if the distributions of the arms are only known to be bounded
(with a known upper bound) and finitely supported. For general bounded
arms, the empirical-likelihood-based upper confidence bounds, which are
the core of the algorithm, still have an adequate level, but obtaining
explicit finite-time regret bounds for the algorithm itself and/or
reducing its computational complexity is still the object of further
investigations; see the discussion in Section~\ref{secKLUCBgeneral}.
The simulation results show that empirical \texttt{KL-UCB} is efficient in
general cases when the distributions are far from being members of
simple parametric families.

In a nutshell, empirical \texttt{KL-UCB} is to be preferred when the
distributions of the arms are not known to belong (or be close) to a
simple parametric family and when the \texttt{kl-UCB}
algorithm is know not to get satisfactory performance---that is, for
instance, when the variance of a $[0,1]$-valued arm with expectation
$\mu$ is much smaller than $\mu(1-\mu)$.


\begin{supplement}
\stitle{Technical proofs}
\slink[doi]{10.1214/13-AOS1119SUPP} 
\sdatatype{.pdf}
\sfilename{aos1119\_supp.pdf}
\sdescription{The supplemental article contains the proofs of the
results stated in the paper.}
\end{supplement}


\printaddresses

\end{document}